
\newcount\secno
\newcount\prmno
\def\section#1{\vskip1truecm
               \global\def\currenvir{section}
               \global\advance\secno by1\global\prmno=0
               {\bf \number\secno. {#1}}
               \smallskip}

\def\subsection{\global\def\currenvir{subsection}
                \global\advance\prmno by1
               \smallskip  \ind{ (\number\secno.\number\prmno) }}
\def\subsec{\global\def\currenvir{subsection}
                \global\advance\prmno by1\smallskip
                { (\number\secno.\number\prmno)\ }}

\def\proclaim#1{\global\advance\prmno by 1
                {\bf #1 \the\secno.\the\prmno$.-$ }}

\long\def\th#1 \enonce#2\endth{%
   \medbreak\proclaim{#1}{\it #2}\global\def\currenvir{th}\smallskip}

\def\rem#1{\global\advance\prmno by 1
{\it #1} \the\secno.\the\prmno$.-$ }

\magnification 1250 \pretolerance=500 \tolerance=1000
\brokenpenalty=5000 \mathcode`A="7041 \mathcode`B="7042
\mathcode`C="7043 \mathcode`D="7044 \mathcode`E="7045
\mathcode`F="7046 \mathcode`G="7047 \mathcode`H="7048
\mathcode`I="7049 \mathcode`J="704A \mathcode`K="704B
\mathcode`L="704C \mathcode`M="704D \mathcode`N="704E
\mathcode`O="704F \mathcode`P="7050 \mathcode`Q="7051
\mathcode`R="7052 \mathcode`S="7053 \mathcode`T="7054
\mathcode`U="7055 \mathcode`V="7056 \mathcode`W="7057
\mathcode`X="7058 \mathcode`Y="7059 \mathcode`Z="705A
\def\spacedmath#1{\def\packedmath##1${\bgroup\mathsurround =0pt##1\egroup$}
\mathsurround#1
\everymath={\packedmath}\everydisplay={\mathsurround=0pt}}
 \spacedmath{2pt}

\def\iso{\vbox{\hbox to .8cm{\hfill{$\scriptstyle\sim$}\hfill}
\nointerlineskip\hbox to .8cm{{\hfill$\longrightarrow $\hfill}} }}
\def\sdir_#1^#2{\mathrel{\mathop{\kern0pt\oplus}\limits_{#1}^{#2}}}

\font\eightrm=cmr8 \font\sixrm=cmr6

\def\pc#1{\tenrm#1\sevenrm}
\def\tx{\kern-1.5pt -}
\def\cqfd{\kern 2truemm\unskip\penalty 500\vrule height 4pt depth 0pt width
4pt\medbreak} 
\def\no{n\up{o}\kern 2pt}
\def\ind{\par\hskip 1truecm\relax}

\font\pal=cmsy7

\def\sp#1{{\cal S}\kern-1pt\raise-1pt\hbox{\pal P}^{}_C(#1)}

\frenchspacing
\input xy
\xyoption{all}
\input amssym.def
\input amssym
\vsize = 25truecm \hsize = 16.1truecm \voffset = -.5truecm
\parindent=0cm
\baselineskip15pt \overfullrule=0pt

\vglue 2.5truecm \font\Bbb=msbm10

\centerline{{\bf On the residue fields of Henselian valued stable 
fields, II}}
\bigskip

\centerline{I.D. Chipchakov\footnote{$^{\ast}$}{Partially
supported by Grant MI-1503/2005 of the Bulgarian Foundation for
Scientific Research.}}
\par
\medskip

\par
\vskip0.75truecm
\centerline{{\bf Introduction} }
\par
\medskip
This paper is a continuation of [Ch2]. Let $E$ be a field, $E _{\rm 
sep}$ a separable closure of $E$, $E ^{\ast }$ the multiplicative 
group of $E$, $d(E)$ the class of finite-dimensional central 
division $E$-algebras, Br$(E)$ the Brauer group of $E$, Gal$(E)$ 
the set of finite Galois extensions of $E$ in $E _{\rm sep}$, and 
NG$(E)$ the set of norm groups $N(M/E)\colon \ M \in {\rm 
Gal}(E)$. Let also $\overline P$ be the set of prime numbers, 
Br$(E) _{p}$ the $p$-component of Br$(E)$, for each $p \in 
\overline P$, $P(E)$ the set of those $p \in \overline P$, for 
which $E$ is properly included in its maximal $p$-extension $E (p)$ 
(in $E _{\rm sep}$), and $\Pi (E)$ the set of all $p ^{\prime } \in 
\overline P$, for which the absolute Galois group $G _{E} = G(E 
_{\rm sep}/E)$ is of nonzero cohomological $p ^{\prime }$-dimension 
cd$_{p'} (G _{E})$. Recall that $E$ is $p$-quasilocal, for a given 
$p \in \overline P$, if $p \not\in P(E)$ or the relative Brauer 
group Br$(F/E)$ equals the subgroup $_{p} {\rm Br}(E) = \{b _{p} 
\in {\rm Br}(E)\colon \ pb _{p} = 0\}$, for every cyclic extension 
$F$ of $E$ of degree $p$. The field $E$ is said to be primarily 
quasilocal (PQL), if it is $p$-quasilocal, for all $p \in \overline 
P$; $E$ is called quasilocal, if its finite extensions are PQL. We 
say that $E$ is a strictly PQL-field, if it is PQL and Br$(E) _{p} 
\neq \{0\}$, for every $p \in P(E)$ (see [Ch4, Corollary 3.7], for 
a characterization of the SQL-property). The purpose of this note 
is to complement the main result of [Ch2] and to shed light on the 
structure of Br$(K)$, where $K$ is an absolutely stable field (in 
the sense of E. Brussel) with a Henselian valuation whose value 
group is totally indivisible. 
\par
\vskip0.6truecm
\centerline{\bf 1. Subfields of central division algebras over 
PQL-fields}
\par
\medskip
Assume that $E$ is a PQL-field and $M \in {\rm Gal}(E)$, such that 
$G(M/E)$ is nilpotent, and let $R$ be an intermediate field of $M/E$.
It follows from Galois theory, [Ch2, Theorem 4.1 (iii)], the 
Burnside-Wielandt characterization of nilpotent finite groups, and 
the general properties of central division algebras and their Schur 
indices (see [KM, Ch. 6, Sect. 2] and [P, Sects. 13.4 and 14.4]) 
that $R$ embeds as an $E$-subalgebra in each algebra $D \in d(E)$ 
of index divisible by $[R\colon E]$. In this Section, we demonstrate 
the optimality of this result in the class of strictly PQL-fields.
\par
\medskip
{\bf Theorem 1.1.} {\it For each nonnilpotent finite group $G$, 
there exists a strictly {\rm PQL}-field $E = E(G)$ with {\rm 
Br}$(E) \cong \hbox{\Bbb Q}/\hbox{\Bbb Z}$ and {\rm Gal}$(E)$ 
containing an element $M$ such that $G(M/E) \cong G$ and $M$ does 
not embed as an $E$-subalgebra in any $\Delta \in d(E)$ of index 
equal to $[M\colon E]$.}
\par
\medskip
{\it Proof.} Let $P(G)$ be the set of prime divisors of the order of
$G$. Our argument relies on the existence (cf. [Ch1, Sect. 4]) of an
algebraic number field $E _{0}$, such that there is $M _{0} \in {\rm 
Gal}(E _{0})$ with $G(M _{0}/E _{0}) \cong G$, and $E _{0}$ possesses 
a system of (real-valued) valuations $\{w(p)\colon \ p \in P(G)\}$ 
satisfying the following conditions:
\par
\medskip
(1.1) (i) $w(p)$ extends the normalized $p$-adic valuation (in the
sense of [CF, Ch. VII]) of the field $\hbox{\Bbb Q}$ of rational
numbers, for each $p \in P(G)$;
\par
(ii) The completion $M _{0,w(p)'}$ lies in Gal$(E _{0,w(p)})$ and 
$G(M  _{0,w(p)'}/E _{0,w(p)})$ is isomorphic to a Sylow 
$p$-subgroup of $G(M _{0}/E _{0})$, whenever $p \in P(G)$ and $w(p) 
^{\prime }$ extends $w(p)$ on $M _{0}$.
\par
\medskip
Let $A _{p}$ be the maximal abelian $p$-extension of $E _{0}$ in $M
_{0}$, for each $p \in P(G)$. It follows from the nonnilpotency of
$G$, the Burnside-Wielandt theorem and Galois theory that $M _{0}/E
_{0}$ has an intermediate field $F \neq E _{0}$ such that $F
\cap A _{p} = E$, for all $p \in P(G)$ (see [Ch1, Sect. 4]). Fix a
prime divisor $\pi $ of $[F\colon E _{0}]$ as well as Sylow $\pi
$-subgroups $H _{\pi } \in {\rm Syl} _{\pi } G(M _{0}/F)$ and $G
_{\pi } \in {\rm Syl} _{\pi } G(M _{0}/E _{0})$ with $H _{\pi }
\subset G _{\pi }$, and denote by $F _{\pi }$ and $E _{\pi }$ the
 maximal extensions of $E _{0}$ in $M _{0}$ fixed by $H _{\pi }$
and $G _{\pi }$, respectively. It is clear from (1.1) (ii) and [CF,
Ch. II, Theorem 10.2] that $E _{\pi }$ has a valuation $t(\pi )$
extending $w(\pi )$ so that $M _{0} \otimes _{E _{\pi }} E _{\pi
, t(\pi )}$ is a field. This means that $t(\pi )$ is uniquely
(up-to an equivalence) extendable to a valuation $\mu (\pi )$ of
$M _{0}$, and allows us to identify $E _{0,w(\pi )}$ with $E _{\pi
,t(\pi )}$. Using repeatedly the Grunwald-Wang theorem [W] and the
normality of maximal subgroups of finite $\pi $-groups (cf. [L, Ch.
I, Sect. 6]), one also obtains that there is a finite extension $K$
of $E _{0}$ in $E _{0} (\pi )$, such that $K \otimes _{E _{0}} E
_{0,w(\pi )}$ is a field isomorphic to $F _{\pi ,\nu (\pi )}$ as an
$E _{0}$-algebra, where $\nu (\pi )$ is the valuation of $F _{\pi
}$ induced by $\mu (\pi )$. These observations indicate that $K 
\cap M _{0} = E _{0}$, i.e. the compositum $M _{0}K$ lies in 
Gal$(K)$ and $G((M _{0}K)/K) \cong G(M _{0}/E _{0})$. At the same 
time, our argument proves that $F$ is dense in $F _{\pi ,\nu (p)}$ 
with respect to the topology induced by $\nu (p)$, $K$ has a unique 
valuation $\kappa (\pi )$ extending $w(\pi )$, and $(M _{0}K) 
_{t(\pi )'} \in {\rm Gal}(K _{\kappa (\pi )})$ with $G((M _{0}K) 
_{t(\pi )'}/K _{\kappa (\pi )}) \cong H _{\pi }$, for any 
prolongation $t(\pi ) ^{\prime }$ of $t(\pi )$ on $M _{0}K$. Let 
now $M(K) = \{\kappa (p)\colon \ p \in \overline P\}$ be a system 
of valuations of $K$ complementing $\kappa (\pi )$ so that $\kappa 
(p)$ extends $w(p)$ or the normalized $p$-adic valuation of 
$\hbox{\Bbb Q}$, depending on whether or not $p \in 
(P(G) \setminus \{\pi \})$. Applying [Ch3, Theorem 2.2] to $K$, 
$\overline P$ and $M(K)$, and using [Ch3, Lemmas 3.1 and 3.2], 
one proves the existence of an extension $E$ of $K$ in $E _{0,{\rm 
sep}}$, such that $P(E) = \overline P$, $M _{0}E := M$ lies in 
Gal$(E)$, $G(M/E) \cong G$, and $E$ possesses a system 
$\{v(p)\colon \ p \in \overline P\}$ of valuations with the 
following properties:
\par
\medskip
(1.2) For each $p \in \overline P$, $v(p)$ extends $\kappa (p)$, $E 
_{v(p)}$ is a PQL-field, $E _{v(p)} (p)$ is $E$-isomorphic to $E 
(p) \otimes _{E} E _{v(p)}$, Br$(E _{v(p)}) _{p} \neq \{0\}$ and 
Br$(E _{v(p')}) _{p} = \{0\}$, for every $p ^{\prime } \in (\overline 
P \setminus \{p\})$. Moreover, if $p \in P(G)$ and $v(p) ^{\prime 
}$ is a valuation of $M$ extending $v(p)$, then $M _{v(p)'} \in 
{\rm Gal}(E _{v(p)})$, $G(M _{v(p)'}/E _{v(p)}) \cong G(M 
_{0,w(p)'}/E _{0,w(p)})\colon \ p \neq \pi $, and $G(M _{v(\pi )'}/E 
_{v(\pi )}) \cong H _{\pi }$.
\par
\medskip \noindent
Hence, by [Ch3, Theorem 2.1], $E$ is a nonreal strictly PQL-field. 
As $E _{v(\pi )}$ is PQL, it follows from [Ch2, Theorem 4.1], 
statements (1.2) and the Brauer-Hasse-Noether and Albert theorem (in 
the form of [Ch3, Proposition 1.2]) that $M$ splits an algebra 
$\Delta \in d(E)$ of $\pi $-primary dimension if and only if 
ind$(\Delta )$ divides the order of $H _{\pi }$. In view of the 
general theory of simple algebras (see [P, Sects. 13.4 and 14.4]), 
this proves Theorem 1.1.
\par
\vskip0.6truecm 
\centerline{\bf 2. Divisible and reduced parts of the 
multiplicative group of an SQL-field}
\par
\medskip
A field $E$ is said to be strictly quasilocal (SQL), if its finite 
extensions are strictly PQL. When this holds and $E$ is almost 
perfect (in the sense of [Ch2, I, (1.8)]), this Section gives a 
Galois-theoretic characterization of the maximal divisible 
subgroup $D(E)$ of $E ^{\ast }$. 
\par
\medskip
{\bf Theorem 2.1.} {\it Let $E$ be an almost perfect {\rm 
SQL}-field and $N _{1} (E)$ the intersection of the groups from 
{\rm NG}$(E)$. Then $N _{1} (E) = D(E ^{\ast })$.}
\par
\medskip
To prove Theorem 2.1 we need the following lemma.
\par
\medskip
{\bf Lemma 2.2.} {\it Let $E$ and $M$ be fields, $M \in{\rm 
Gal}(E)$, and let $P(M/E)$ be the set of prime divisors of 
$[M\colon E]$. Then $N(M/E) \subseteq N(M/F)$, for every 
intermediate field $F$ of $M/E$. Moreover, if $E _{p}$ is the fixed 
field of some Sylow $p$-subgroup of $G(M/E)$, then $N(M/E) = \cap _{p 
\in P(M/E)} N(M/E _{p})$.}
\par
\medskip
{\it Proof.} It is easily verified that if $[F\colon E] = m$, 
$\{\tau _{j}\colon \ j = 1, \dots , m\}$ is the set of 
$E$-embeddings of $F$ into $M$, and $\sigma _{u}$ is an 
automorphism of $M$ extending $\tau _{u}$, for each index $j$, then 
$G(M/E) = \cup _{j=1} ^{m} G(M/F)\sigma _{j} ^{-1}$. This yields $N 
_{E} ^{M} (\gamma ) = N _{F} ^{M} (\prod _{j=1} ^{m} \sigma _{j} 
^{-1} (\gamma ))$, for any $\gamma \in M ^{\ast }$, which proves 
that $N(M/F) \subseteq N(M/E)$. We show that $\cap _{p \in P(M/E)}$ 
\par \noindent
$N(M/E _{p}) := N _{0} (M/E) \subseteq N(M/E)$. Take an element 
$\alpha \in N _{0} (M/E)$ and put $[E _{p}\colon E]$
\par \noindent
$= m _{p}$, for each $p \in P(M/E)$. It is known that $p \not\vert 
m _{p}$, for any $p \in P(M/E)$, which implies consecutively that 
${\rm g.c.d.}(m _{p}\colon p \in P(M/E)) = 1$, $\alpha \in E ^{\ast 
}$ and $\alpha ^{m _{p}} \in N(M/E)$, for every $p \in P(M/E)$. 
Thus it turns out that $\alpha \in N(M/E)$, so Lemma 2.2 is proved.
\par
\medskip
{\it Proof of Theorem 2.1.} The inclusion $D(E ^{\ast }) \subseteq N 
_{1} (E)$ is obvious, so our objective is to prove the converse. We 
show that $N _{1} (E) \subseteq N _{1} (E) ^{p ^{n}}$, for every $p 
\in \overline P$ and $n \in \hbox{\Bbb N}$. It is clearly sufficient 
to consider the special case where $p \in \Pi (E)$. Then there is a 
finite extension $F _{p}$ of $E$ in $E _{\rm sep}$, such that Br$(F 
_{p}) _{p} \neq \{0\}$ and $p \not\vert [F _{p}\colon E]$. Note 
also that each finite extension of $F _{p}$ in $E ^{\ast }$ is 
included in a field $\widetilde F _{p} \in {\rm Gal}(E)$, so it 
follows from Lemma 2.2 that $N _{1} (E) \subseteq N _{1} (F _{p})$. 
At the same time, the norm mapping $N _{E} ^{F _{p}}$ clearly
induces homomorphisms $D(F _{p} ^{\ast }) \to D(E ^{\ast })$ and $N 
_{1} (F _{p}) \to N _{p} (E)$ (for the latter, see [L, Ch. VIII, 
Sect. 5]) as well as an automorphism of the maximal $p$-divisible 
subgroup of $E ^{\ast }$. These observations allow one to assume 
additionally that $p \in \overline P(E)$ (and Br$(E) _{p} \neq 
\{0\}$). If $p = {\rm char}(E)$, our assertion is implied by [Ch2, 
Lemma 8.4], so we turn to the case of $p \neq {\rm char}(E)$. Let 
$\varepsilon _{p}$ be a primitive $p$-th root of unity in $E _{\rm 
sep}$. Then $[E(\varepsilon _{p})\colon E]$ divides $p - 1$ (cf. 
[L, Ch. VIII, Sect. 3]), whence our considerations reduce further 
to the case in which $\varepsilon _{p} \in E$. Now the proof of 
Theorem 2.1 is completed by applying the following two lemmas. 
Before stating them, let us recall that the character group 
$C(Y/E)$ of $G(Y/E)$ is an abelian torsion group, for every Galois 
extension $Y/E$ (see [K, Ch. 7, Sect. 5]). Therefore, by [F, 
Theorem 24.5], the divisible part $D(Y/E)$ of $C(Y/E)$ is a direct 
summand in $C(Y/E)$, i.e. $C(Y/E)$ is isomorphic to the direct sum 
$D(Y/E) \oplus R(Y/E)$, where $R(Y/E) \cong C(Y/E)/D(Y/E)$ is a 
maximal reduced subgroup of $C(Y/E)$. 
\par
\medskip
{\bf Lemma 2.3.} {\it Let $E$ be a $p$-quasilocal field containing 
a primitive $p$-th root of unity $\varepsilon $, for some prime $p 
\in P(E)$. Assume also that $r _{p} (E)$ is the group of roots of 
unity in $E$ of $p$-primary degrees. Then:}
\par
(i) $R(E (p)/E) = \{0\}$ {\it if and only if {\rm Br}$(E) _{p} = 
\{0\}$ or $r _{p} (E)$ is infinite;}
\par
(ii) {\it If {\rm Br}$(E) _{p} \neq \{0\}$, $d$ is the dimension of 
$_{p} {\rm Br}(E)$ as a vector space over the field $\hbox{\Bbb F} 
_{p}$ with $p$ elements, and $R _{p} (E)$ is of order $p ^{\mu }$, for 
some $\mu \in \hbox{\Bbb N}$, then $D(E (p)/E) = p ^{\mu }C(E (p)/E)$ 
and $R(E (p)/E)$ is presentable as a direct sum of cyclic groups of 
order $p ^{\mu }$, indexed by a set of cardinality $d$.}
\par
\medskip
{\it Proof.} Statement (i) is a well-known consequence of Kummer 
theory, the Merkurjev-Suslin theorem and Galois cohomology (see 
[MS, (11.5), (16.1)] and [S, Ch. I, 3.4 and 4.2]), so we assume 
further that $d > 0$ and $r _{p} (E)$ is of order $p ^{\mu }$, for 
some $\mu \in \hbox{\Bbb N}$. We first show that $p ^{\mu }C(E 
(p)/E)$ is divisible. Fix a primitive $p ^{\mu }$-th root of unity 
$\varepsilon _{\mu }$ in $E$ and a subset $S = \{\varepsilon 
_{n}\colon \ n \in \hbox{\Bbb N}, n \ge \mu + 1\}$ of $E (p)$ so 
that $\varepsilon _{n+1} ^{p} = \varepsilon _{n}$, for each index 
$n$. Also, let $E _{\infty } = E(S)$ and $E _{n} = E(\varepsilon 
_{n})$, for every integer $n > \mu $. Suppose first that $p > 2$ or 
$\mu \ge 2$. Then $E _{\infty }/E$ is a $\hbox{\Bbb Z} 
_{p}$-extension and $E _{n}$ is the unique subextension of $E$ in 
$E _{\infty }$ of degree $p ^{n-\mu }$. Since $\hbox{\Bbb Z} _{p}$ 
is a projective profinite group (see [G, Theorem 1]), this enables 
one to deduce from Galois theory that $E (p)$ possesses a subfield 
$E ^{\prime }$ such that $E ^{\prime } \cap E _{\infty } = E$ and 
$E _{\infty }E ^{\prime } = E (p)$. In addition, it becomes clear 
that $C(E (p)/E)$ is isomorphic to the direct sum $\hbox{\Bbb Z} (p 
^{\infty }) \oplus C(F/E)$, where $\hbox{\Bbb Z} (p ^{\infty })$ is 
the quasicyclic $p$-group (identified with the character group of 
$\hbox{\Bbb Z} _{p}$) and $F$ is the maximal abelian extension of 
$E$ in $E ^{\prime }$. Let now $\Phi $ be a cyclic extension of $E$ 
in $F$. As $E$ is $p$-quasilocal, an element $\beta _{n} \in E _{n} 
^{\ast }$ lies in $N((\Phi E _{n})/E _{n})$, for an arbitrary 
integer $n \ge \mu $, if and only if the norm $N _{E} ^{E _{n}} 
(\beta _{n}) \in N(\Phi /E)$. In particular, this proves that 
$\varepsilon _{\mu } \in N(\Phi /E)$ if and only if $\varepsilon 
_{n} \in N((\Phi E _{n})/E _{n})$. In view of the results of [FSS, 
Sect. 2], this means that $p ^{\mu }C(F/E) = p ^{n}C(F/E)$, for 
every index $n \ge \mu $. In other words, $C(F/E)$ is divisible, so 
$p ^{\mu }C(E (p)/E)$ has the same property, as claimed.
\par
Our objective now is to prove the divisibility of $2C(E (2)/E)$,
under the hypothesis that $p = 2$ and $\mu = 1$. If $E _{\infty } =
E(\sqrt{-1}) = E _{2}$, this is contained in [FSS, Proposition 2],
so we assume further that $E _{\infty } \neq E _{2}$, i.e. $E _{2}
= E _{\nu } \neq E _{\nu +1}$, for some integer $\nu \ge 2$. Let
$R$ be a cyclic extension of $E$ in $E (2)$. By Albert's theorem 
(see [FSS, Sect. 2]), we have $C(R/E) \subset 2C(E (2)/E)$ if and 
only if $-1 \in N(R/E)$. We show that if $C(R/E) \subset 2C(E 
(2)/E)$, then $\varepsilon _{\nu } \in N((RE _{2})/E _{2})$. Since 
$N _{E} ^{E _{2}} (\varepsilon _{\nu })$ equals  $1$ or $-1$, this 
follows from [Ch2, Lemma 4.2 (iii)] in case $R \cap E _{2} = E$. 
Suppose further that $E _{2} \subseteq R$ and fix a generator 
$\sigma $ of $G(F/E)$. Applying [Ch2, Theorem 4.1] one obtains that 
the cyclic $E _{2}$-algebra $(R/E _{2}, \sigma ^{2}, \varepsilon 
_{\nu })$ is similar to $(R/E, \sigma , c) \otimes _{E} E _{2}$, 
for some $c \in E ^{\ast }$. In view of [P, Sect. 14.7, Proposition 
b], this ensures that $(R/E _{2}, \sigma ^{2}, \varepsilon _{\nu }) 
\cong (R/E _{2}, \sigma ^{2}, c)$ as $E _{2}$-algebras. Therefore, 
[P, Sect. 15.1, Proposition b] yields $c\varepsilon _{\nu } ^{-1} 
\in N(R/E _{2})$ and $c ^{2}N _{E} ^{E _{2}} (\varepsilon _{\nu }) 
\in N(R/E)$. Since $-1 \in N(R/E)$, this means that $c ^{2} \in 
N(R/E)$ as well. Observing now that the corestriction homomorphism 
of Br$(E _{2})$ into Br$(E)$ maps the similarity class $A _{c}$ of 
$(R/E _{2}, \sigma ^{2}, c)$ into the similarity class of $(R/E, 
\sigma , c ^{2})$ (cf. [T, Theorem 2.5]), one concludes that $A 
_{c}$ lies in the kernel Ker$_{E _{2}/E}$ of this homomorphism. At 
the same time, by [Ch2, Lemma 4.2 (i)], we have Br$(E _{2}) _{2} 
\cap {\rm Ker}_{E _{2}/E} = \{0\}$. In particular, $A _{c} = 0$, so 
it follows from [P, Sect. 15.1, Proposition b] that $\varepsilon 
_{\nu } \in N(R/E _{2})$, as claimed. The obtained result, combined 
with [AFSS, Theorem 3], implies that $2C(E (2)/E) = 2 ^{\nu }C(E 
(2)/E) = 4C(E (2)/E)$, which proves the divisibility of $2C(E 
(2)/E)$. 
\par
In order to complete the proof of Lemma 2.3, it remains to be seen
that the quotient group $C(E (p)/E)/p ^{\mu }C(E (p)/E)$ is
isomorphic to the direct sum of cyclic groups of order $p ^{\mu }$,
indexed by a set $I$ of cardinality $d$ (see [F, Theorem 24.5]).
Since, by Kummer's theory, $p ^{\mu -1}C(E (p)/E)$ includes the 
group $X _{p} (E) = \{\chi \in C(E (p)/E)\colon \ p\chi = 0\}$, it 
is sufficient to show that $(X _{p} (E) + p ^{\mu }C(E (p)/E))/p 
^{\mu }C(E (p)/E) := \overline X _{p} (E)$ is isomorphic to $_{p} 
{\rm Br}(E)$. Applying Albert's theorem and elementary properties 
of symbol $E$-algebras, and taking into account that $\overline X 
_{p} (E) \cong X _{p} (E)/(X _{p} (E) \cap p ^{\mu }C(E (p)/E))$, 
one obtains that $\overline X _{p} (E) \cong E ^{\ast }/N(E _{\mu 
+1}/E)$. At the same time, the cyclicity of $E _{\mu +1}/E$ ensures 
that $E ^{\ast }/N(E _{\mu +1}/E) \cong {\rm Br}(E _{\mu +1}/E)$ 
(cf. [P, Sect. 15.1, Proposition b]), and by the $p$-quasilocal 
property of $E$, we have Br$(E _{\mu +1}/E) = \ _{p} {\rm Br}(E)$, so 
our proof is complete.
\par
\medskip
{\bf Lemma 2.4.} {\it Let $E$ be a $p$-quasilocal field containing 
a primitive $p$-th root of unity $\varepsilon $, and such that 
Br$(E) _{p} \neq \{0\}$, and let $\Omega _{p} (E)$ be the set of 
finite abelian extensions of $E$ in $E (p)$. Then the intersection 
$N _{p} (E)$ of the norm groups $N(M/E)\colon \ M \in \Omega _{p} 
(E)$, coincides with the maximal $p$-divisible subgroup of $E 
^{\ast }$.}
\par
\medskip
{\it Proof.} It is clearly sufficient to show that $N _{p} (E) 
\subseteq N _{1} (E) ^{p ^{n}}$, for each $n \in \hbox{\Bbb N}$. 
Denote by $r _{p} (E)$ the group of roots of unity in $E$ of 
$p$-primary degrees and fix an algebra $D \in d(E)$ of index $p$ 
(the existence of $D$ is guaranteed by the assumption that Br$(E) 
_{p} \neq \{0\}$ and the Merkurjev-Suslin theorem [MS, (16.1)]). As 
$E$ is $p$-quasilocal, one obtains from Kummer theory that, for 
each $c \in E ^{\ast } \setminus E ^{\ast p}$, there is $c ^{\prime 
} \in E ^{\ast } \setminus E ^{\ast p}$, such that $D$ is 
$E$-isomorphic to the symbol $E$-algebra $A _{\varepsilon } (c, c 
^{\prime }; E)$. Hence, $c \not\in N(E _{c'}/E)$, where $E _{c'}$ 
is the extension of $E$ in $E (p)$ obtained by adjunction of a 
$p$-th root of $c ^{\prime }$ (see [P, Sect. 15.1, Proposition 
b]). The obtained result implies $N _{p} (E) \subseteq E ^{\ast 
p}$. Arguing in a similar manner and applying general properties of 
cyclic $E$-algebras (cf. [P, Sect. 15.1, Corollary b]), one deduces 
that $N _{p} (E) \subseteq E ^{\ast p ^{m}}$, if $E$ contains a 
primitive $p ^{m}$-th root of unity (and proves the lemma in case 
$r _{p} (E)$ is infinite). Suppose further that $r _{p} (E)$ is of 
order $p ^{\mu }$, for some $\mu \in \hbox{\Bbb N}$, and take $C(E 
(p)/E)$, $D(E (p)/E)$ and $R(E (p)/E)$ as in Lemma 2.3. Also, let 
$c \in (E ^{\ast p ^{m}} \cap N _{p} (E))$, for some integer $m \ge 
\mu $. We prove Lemma 2.4 by showing that $c \in (E ^{\ast p 
^{m+\mu }} \cap N _{p} (E) ^{p ^{m}})$. Fix a primitive $p ^{\mu 
}$-th root of unity $\delta _{\mu } \in E$, take an element $c _{m} 
\in E ^{\ast }$ so that $c _{m} ^{p ^{m}} = c$, and for any 
$\lambda \in E ^{\ast }$, denote by $E _{\lambda }$ the extension 
of $E$ in $E (p)$ obtained by adjunction of a $p ^{\mu }$-th root 
of $\lambda $. Using again [P, Sect. 15.1, Corollary b], one 
obtains that $c _{m} \in N(F/E)$ whenever $F$ is an intermediate 
field of a $\hbox{\Bbb Z} _{p}$-extension of $E$. Hence, by 
Albert's theorem (see [FSS, Sect. 2]), Kummer theory and the basic 
properties of symbol $E$-algebras of dimension $p ^{2\mu }$, $N(E 
_{\delta _{\mu }})/E) \subseteq N(E _{c _{m}})/E)$. Since Br$(E) 
_{p} \neq \{0\}$, $E$ admits (one-dimensional) local $p$-class 
field theory [Ch4, Theorem 3.1], which means that $E _{c _{m}} 
\subseteq E _{\delta _{\mu }}$. In view of Kummer theory, the 
obtained result yields $c _{m} = \delta _{\mu } ^{k}t _{m} ^{p 
^{\mu }}$, for some $k \in \hbox{\Bbb N}$, $t _{m} \in E ^{\ast }$. 
Therefore, we have $t _{m} ^{p ^{m+\mu }} = (t _{m} ^{p ^{\mu }}) 
^{p ^{m}} = c$. Let now $M$ be an arbitrary finite abelian 
extension of $E$ in $E (p)$. Then it follows from Galois theory and 
Lemma 2.3 that $M$ is a subfield of the compositum of finitely many 
cyclic extensions $D _{i}, i \in I$, and $R _{j}, j \in J$, of $E$ 
in $E (p)$, such that $[R _{j}\colon E] = p ^{\mu }$, for every $j 
\in J$, and $D _{i} \subseteq \Delta _{i}$, where $\Delta _{i}$ is 
a $\hbox{\Bbb Z} _{p}$-extension of $E$ in $E (p)$, for each index 
$i$. Hence, by [Ch4, Theorem 3.1], $c _{m} ^{p ^{\mu }} \in 
N(\widetilde M/E) \subseteq N(M/E)$, i.e. $c _{m} ^{p ^{\mu }} \in 
N _{p} (E)$. Replacing $c _{m}$ by $t _{m}$ and arguing as above, 
one obtains that $t _{m} ^{p ^{\mu }} \in N _{p} (E)$ and $c \in N 
_{p} (E) ^{p ^{m}}$. Thus the assertion that $c \in (E ^{\ast p 
^{m+\mu }} \cap N _{p} (E) ^{p ^{m}})$ is proved. It is now easy to 
see that $N _{p} (E) = N _{p} (E) ^{p ^{n}}$, for any $n \in 
\hbox{\Bbb N}$, as required by Lemma 2.4. 
\par
\medskip
{\bf Remark 2.5.} Analyzing the proof of Theorem 2.1, one obtains 
that its conclusion remains valid, if $E$ is a quasilocal perfect 
field and cd$_{p} (G _{E}) \neq 1$, for all $p \in (\overline P 
\setminus {\rm char}(E))$. One also proves that $R(E ^{\ast }) 
\cong \Phi (E) \oplus R _{0} (E)$, where $\Phi (E)$ is a 
torsion-free group and $R _{0} (E)$ is the direct sum of the groups 
$r _{p} (E)$, indexed by those $p \in \Pi (E)$, for which $r _{p} 
(E)$ is nontrivial and finite. Moreover, it turns out that if 
Br$(E) _{p} \neq \{0\}$, for every $p \in (\Pi (E) \setminus {\rm 
char}(E))$, then $D(E ^{\ast })$ equals the intersection of the 
groups $N(M/E)$, defined over all $M \in {\rm Gal}(E)$, such that 
$G(M/E)$ lies in any fixed class $\chi $ of finite groups, which is 
closed under the formation of subgroups, homomorphic images and 
finite direct products, and which contains all finite metabelian 
groups of orders not divisible by any $p \in (\overline P \setminus 
\Pi (E))$.  
\par
\vskip0.6truecm
\centerline{\bf 3. The reduced part of the Brauer group of an 
equicharacteristic Henselian}
\par
\centerline{\bf valued absolutely stable field with a totally 
indivisible value group}
\par
\medskip
In what follows, our notation agrees with that of Lemma 2.3 and its 
proof. For each Henselian valued field $(K, v)$, $\widehat K$ and 
$v(K)$ denote the residue field and the value group of $(K, v)$, 
respectively. In this Section we announce the following 
characterization of the reduced components of the Brauer groups of 
the fields pointed out in its title: 
\par
\medskip
{\bf Theorem 3.1.} {\it An abelian torsion group $T$ is isomorphic 
to the reduced part of {\rm Br}$(K)$, for some absolutely stable 
field $K = K(T)$ with a Henselian valuation $v$ such that {\rm 
char}$(K) = {\rm char}(\widehat K)$ and the value group $v(K)$ is 
totally indivisible (i.e. $v(K) \neq pv(K)$, for every $p \in 
\overline P$), if and only if the $p$-component $T _{p}$ of $T$ is 
presentable as a direct sum of cyclic groups of one and the same 
order $p ^{n _{p}}$, for each $p \in \overline P$.} 
\par
\medskip
The rest of this Section is devoted to the proof of the necessity 
in Theorem 3.1 (the sufficiency will be proved elsewhere in a more 
general situation covering the case where char$(K) \neq {\rm 
char}(\widehat K)$, $v$ is Henselian and discrete, and $\widehat K$ 
is perfect). We begin with a description of some known basic 
relations between Br$(K) _{p}$, Br$(\widehat K) _{p}$ and 
$C(\widehat K (p)/\widehat K)$ (proved for convenience of the 
reader).
\par
\medskip
{\bf Lemma 3.2.} {\it Let $(K, v)$ be a Henselian valued field
with $v(K) \neq pv(K)$, for some $p \in \overline P$, and let $\pi 
$ be an element of $K ^{\ast }$ of value $v(\pi ) \not\in pv(K)$. 
For each $\chi \in C(\widehat K (p)/\widehat K)$, denote by 
$\widetilde L _{\chi }$ the extension of $\widehat K$ in $\widehat 
K _{\rm sep}$ corresponding by Galois theory to the kernel {\rm 
Ker}$(\chi )$, and by $\tilde \sigma _{\chi }$ the generator of 
$G(\widetilde L _{\chi }/\widehat K)$ satisfying the equality $\chi 
(\tilde \sigma _{\chi }) = (1/[L _{\chi }\colon \widehat K]) + 
\hbox{\Bbb Q}/\hbox{\Bbb Z}$. Assume also that $L {\chi }$ is the 
inertial lift of $\widetilde L _{\chi }$ (over $\widehat K$) in $K 
_{\rm sep}$, $\eta _{\chi }$ is the canonical isomorphism of $G(L 
_{\chi }/K)$ on $G(\widetilde L _{\chi }/\widehat K)$, and $\sigma 
_{\chi }$ is the preimage of $\tilde \sigma _{\chi }$ in $G(L 
_{\chi }/K)$ with respect to $\eta _{\chi }$. Then the mapping $W 
_{\pi }$ of $C(\widehat K (p)/\widehat K)$ into {\rm Br}$(K) _{p}$ 
defined by the rule $W _{\pi } (\chi ) = [(L _{\chi }/K, \sigma 
_{\chi }, \pi )]\colon \ \chi \in C(\widehat K (p)/\widehat K)$, is 
an injective group homomorphism.}
\par
\medskip
{\it Proof.} It is clearly sufficient to show that $[(L _{\chi
_{1} + \chi _{2}}/K, \sigma _{\chi _{1} + \chi _{2}}, \pi )] = [(L
_{\chi _{1}}/L, \sigma _{\chi _{1}}, \pi )]$
\par \noindent
$+ [(L _{\chi _{2}}/K, \sigma _{\chi _{2}}, \pi )]$, where $\chi 
_{1}, \chi _{2} \in C(\widehat K (p)/\widehat K)$, in each of the 
following special cases:
\par
\medskip
(3.1) (i) $\chi _{2} \in \langle \chi _{1}\rangle $, i.e. $L
_{\chi _{2}} \subseteq L _{\chi _{1}}$;
\par
(ii) $\langle \chi _{1}\rangle \cap \langle \chi _{2}\rangle =
\{0\}$, i.e. $L _{\chi _{1}} \cap L _{\chi _{2}} = K$.
\par
\medskip
In case (3.1) (i), this follows from the general properties of
cyclic algebras (cf. [P, Sect. 15.1, Corollary b and Proposition
b]). Assuming that $L _{\chi _{1}} \cap L _{\chi _{2}} = K$ and
$[L _{\chi _{1}}\colon K] \ge [L _{\chi _{2}}\colon K]$, one
obtains that the $K$-algebras $(L _{\chi _{1}}/K, \sigma _{\chi
_{2}}, \pi ) \otimes _{K} (L _{\chi _{2}}/K, \sigma _{\chi
_{2}},$
\par \noindent
$\pi )$ and $(L _{\chi _{1} + \chi _{2}}/K, \sigma _{\chi _{1} + 
\chi _{2}}, \pi ) \otimes _{K} (L _{\chi _{2}}/K, \sigma _{\chi 
_{2}}, 1)$ are isomorphic, which completes our proof.
\par
\medskip
With assumptions being as in the lemma, let $v(K) \neq pv(K)$, for 
some prime $p \neq {\rm char}(\widehat K)$, and let $B _{p}$ be a 
basis and $n(p)$ the dimension of $v(K)/pv(K)$ as a vector space 
over $\hbox{\Bbb F} _{p}$. Denote by $C(\widehat K (p)/\widehat K) 
^{n(p)}$ the direct sum of some isomorphic copies of $C(\widehat K 
(p)/\widehat K)$, indexed by $B _{p}$. Fix a linear ordering $\le $ 
on $B _{p}$ and put $J _{p} = \{(c _{p}, d _{p}) \in (B _{p} \times 
B _{p})\colon \ c _{p} < d _{p}\}$, and in case $n(p) \ge 2$, take 
a direct sum $R _{p} (\widehat K) ^{J _{p}}$ of isomorphic copies 
of $R _{p} (\widehat K)$, indexed by $J _{p}$. Then the 
Ostrowski-Draxl theorem, the Jacob-Wadsworth decomposition 
lemmas (see [JW]) and Lemma 3.2 imply the following variant of the 
Scharlau-Witt theorem [Sch]:
\par
\medskip
(3.2) (i) Br$(K) _{p} \cong {\rm Br}(\widehat K) _{p} \oplus 
C(\widehat K (p)/\widehat K) ^{n(p)}$ unless $n(p) \ge 2$ and $R 
_{p} (\widehat K) \neq \{1\}$;
\par
(ii) Br$(K) _{p} \cong {\rm Br}(\widehat K)
_{p} \oplus C(\widehat K (p)/\widehat K) ^{n(p)} \oplus R _{p}
(\widehat K) ^{J _{p}}$, if $n(p) \ge 2$ and $R _{p} (\widehat K) 
\neq \{1\}$.
\par
\medskip
Suppose now that $K$ is a Henselian valued absolutely stable 
field such that $v(K)$ is totally indivisible. By [Ch2, Proposition 
2.3], $\widehat K$ is quasilocal, so the necessity in Theorem 3.1 
can be deduced from (3.2), the divisibility of Br$(F) _{q}$, for 
any field $F$ of characteristic $q > 0$ (Witt's theorem, see [Dr, 
Sect. 15]), and the following lemma. Before stating it, note that 
the extension of $E$ in $E _{\rm sep}$ obtained by adjunction of a 
primitive $p$-th root of unity, for some $p \in \overline P$, is 
cyclic of degree dividing $p - 1$ (cf. [L, Ch. VIII, Sect. 3]). 
\par
\medskip
{\bf Lemma 3.3.} {\it Let $E$ be a quasilocal field not containing 
a primitive $p$-th root of unity, for a given prime $p \neq {\rm 
char}(E)$, and let $\varepsilon $ be such a root in $E _{\rm sep}$. 
Fix a generator $\varphi $ of $G(E(\varepsilon )/E)$, take an 
integer $s$ so as to satisfy the equality $\varphi (\varepsilon ) = 
\varepsilon ^{s}$, and put $B _{s} = \{b \in \ _{p} {\rm 
Br}(E(\varepsilon ))\colon \ \varphi (b) = sb\}$. Then:}
\par
(i) $R(E (p)/E) = \{0\}$ {\it if and only if $r _{p} (E(\varepsilon 
))$ is infinite or the group $B _{s}$ is trivial;}
\par
(ii) {\it If $B _{s}$ is nontrivial and of dimension $d$ as an 
$\hbox{\Bbb F} _{p}$-vector space, and if $r _{p} (E(\varepsilon 
))$ is of order $p ^{\mu }$, for some $\mu \in \hbox{\Bbb N}$, then 
$D(E (p)/E) = p ^{\mu }C(E (p)/E)$ and $R(E (p)/E)$ is a direct sum 
of cyclic groups of order $p ^{\mu }$, indexed by a set of 
cardinality $d$.}
\par
\medskip
{\it Proof.} Let $\Lambda = \{\lambda \in E(\varepsilon ) ^{\ast 
}\colon \ \varphi (\lambda )\lambda ^{-s} \in E(\varepsilon ) 
^{\ast p}\}$. It is known (Albert, see [A, Ch. IX, Theorem 15]) 
that an extension $L$ of $E$ in $E _{\rm sep}$ is cyclic of degree 
$p$ if and only if $L(\varepsilon )$ is generated over 
$E(\varepsilon )$ by a $p$-th root of an element $\lambda _{0} \in 
(\Lambda \setminus E(\varepsilon ) ^{\ast p})$. Hence, the Kummer 
isomorphism of $X _{p} (E(\varepsilon ))$ on $E(\varepsilon ) 
^{\ast }/E(\varepsilon ) ^{\ast p}$ induces an isomorphism of $X 
_{p} (E)$ on $\Lambda /E(\varepsilon ) ^{\ast p}$. Observe also 
that the similarity class of the symbol $E(\varepsilon )$-algebra 
$A _{\varepsilon } (\lambda _{1}, \lambda _{2}; E(\varepsilon ))$ 
lies in $B _{i}$ whenever $\lambda _{j} \in \Lambda $, $j = 1, 2$. 
Conversely, by the $p$-quasilocal property of $E(\varepsilon )$, 
each element of $B _{s}$ is presented by a symbol $E(\varepsilon 
)$-algebra determined by a pair of elements of $\Lambda $. The 
noted results enable one to prove the lemma arguing as in the proof 
of Lemma 2.3.
\vskip1cm \centerline{ REFERENCES} \vglue15pt\baselineskip12.8pt
\def\num#1{\smallskip\item{\hbox to\parindent{\enskip [#1]\hfill}}}
\parindent=1.38cm
\par
\medskip
\par
\num{A} A.A. {\pc ALBERT}, {\sl Modern Higher Algebra.} Univ. of 
Chicago Press, XIV, Chicago, Ill., 1937.
\par
\num{AFSS} J.K. {\pc ARASON}, B. {\pc FEIN}, M. {\pc SCHCHER}, J. 
{\pc SONN}, {\sl Cyclic extensions of $K(\sqrt{-1})$.} Trans. Amer. 
Math. Soc. {\bf 313}, No 2 (1989), 843-851.
\par
\num{CF} J.W.S. {\pc CASSELS}, A. {\pc ${\rm FR}\ddot o{\rm HLICH}$} 
(Eds.), {\sl Algebraic Number Theory.} Academic Press, London-New 
York, 1967.
\par
\num{Ch1} I.D. {\pc CHIPCHAKOV}, {\sl On nilpotent Galois groups 
and the scope of the norm limitation theorem in one-dimensional 
abstract local class field theory.} In: Proc. of ICTAMI 05, Alba 
Iulia, Romania, 15.9.-18.9, 2005; Acta Univ. Apulensis, No 10 
(2005), 149-167.
\par
\num{Ch2} I.D. {\pc CHIPCHAKOV}, {\sl On the residue fields of Henselian valued 
stable fields.} Preprint, v. 3 (to appear at 
www.math.arXiv.org/math.RA/0412544).
\par
\num{Ch3} I.D. {\pc CHIPCHAKOV}, {\sl Algebraic extensions of global fields
admitting one-dimensional local class field theory.} Preprint
(available at www.arXiv.org/
\par 
math.NT/0504021).
\par
\num{Ch4} I.D. {\pc CHIPCHAKOV}, {\sl One-dimensional abstract local class field 
theory.} Preprint, v. 3 (available at www.arXiv.org/math.RA/0506515).
\par
\num{Dr} P.K. {\pc DRAXL}, {\sl Skew Fields.} London Math. Soc. Lect. Note Series, 
81, Cambridge etc., Cambridge Univ. Press, 1983. 
\par
\num{FSS} B. {\pc FEIN}, D. {\pc SALTMAN}, M. {\pc SCHACHER}, {\sl 
Heights of cyclic field extensions.} Bull. Soc. Math. Belg., Ser. 
A, {\bf 40} (1988), 213-223.
\par
\num{F} L. {\pc FUCHS}, {\sl Infinite Abelian Groups.} Academic Press,  New
York-London, 1970.
\par
\num{G} K.W. {\pc GRUENBERG}, {\sl Projective profinite groups.} J. Lond. Math. 
Soc. {\bf 42} (1967), 155-165.
\par
\num{KM} M.I. {\pc KARGAPOLOV}, Yu.I. {\pc MERZLJAKOV}, {\sl Fundamentals of the
Theory of Groups, 2nd Ed.} Nauka, Moscow, 1977 (Russian: Engl.
Transl. in Graduate Texts in Math. {\bf 62}, Springer-Verlag, New
York-Heidelberg-Berlin, 1979).
\par
\num{K} G. {\pc KARPILOVSKY}, {\sl Topics in Field Theory.} North-Holland Math.
Stud. {\bf 155}, Amsterdam, 1989.
\par
\num{JW} B. {\pc JACOB}, A.R. {\pc WADSWORTH}, {\sl Division 
algebras over Henselian fields.} J. Algebra {\bf 128} (1990), 
528-579.
\par
\num{L} S. {\pc LANG}, {\sl Algebra}, Addison-Wesley, Reading, MA, 1965.
\par
\num{MS} A.S. {\pc MERKURJEV}, A.A. {\pc SUSLIN}, {\sl $K$-cohomology  of
Brauer-Severi varieties and norm residue homomorphisms.} Izv. Akad.
Nauk SSSR {\bf 46} (1982), 1011-1046  (Russian: Engl. transl. in
Math. USSR Izv. {\bf 21} (1983), 307-340).
\par
\num{P} R. {\pc PIERCE}, {\sl Associative Algebras.} Graduate Texts 
in Math. {\bf 88}, Springer-Verlag, New York-Heidelberg-Berlin, 
1982.
\par
\num{Sch} W. {\pc SCARLAU}, {\sl $\ddot U$ber die Brauer-Gruppe eines 
Hensel-K$\ddot o$rpers.} Abh. Math. Semin. Univ. Hamb. {\bf 33} 
(1969), 243-249.
\par
\num{S} J.-P. {\pc SERRE}, {\sl Cohomologie Galoisienne.} Lect. Notes in Math. 
{\bf 5}, Springer-Verlag, Berlin-Heidelberg-New York, 1965. 
\par
\num{T} J.-P. {\pc TIGNOL}, {\sl On the corestriction of central simple algebras.} 
Math. Z. {\bf 194} (1987), 267-274. 
\par
\num{W} S. {\pc WANG}, {\sl On Grunwald's theorem.} Ann. Math. (2), 51 (1950),
471-484.
\vskip1cm
\def\pc#1{\eightrm#1\sixrm}
\hfill\vtop{\eightrm\hbox to 5cm{\hfill Ivan {\pc
CHIPCHAKOV}\hfill}
 \hbox to 5cm{\hfill Institute of Mathematics and Informatics\hfill}\vskip-2pt
 \hbox to 5cm{\hfill Bulgarian Academy of Sciences\hfill}
\hbox to 5cm{\hfill Acad. G. Bonchev Str., bl. 8\hfill} \hbox to
5cm{\hfill 1113 {\pc SOFIA,} Bulgaria\hfill}}
\end